# THE COALESCENT EFFECTIVE SIZE OF AGE-STRUCTURED POPULATIONS[1]

By Serik Sagitov and Peter Jagers

*Chalmers University of Technology*

We establish convergence to the Kingman coalescent for a class of age-structured population models with time-constant population size. Time is discrete with unit called a year. Offspring numbers in a year may depend on mother's age.

**1. Introduction.** The well-known *coalescent* process describes how family lines merge in a sample from a large population, when time is traced backward. It is a continuous time Markov chain which keeps record of branches starting from $n$ leaves and going through $n-1$ pairwise mergers toward the root of a so-called ultrametric tree. The number of branches is reduced from $k, 2 \le k \le n$, to $k-1$ at the rate $\binom{k}{2}$. At each reduction, a random pair of branches is replaced by a single branch.

Initially [9], the coalescent was obtained as an approximation of the genealogical tree of the Wright–Fisher model (WFM) with a large population size $N$ and the unit of coalescent time corresponding to $N$ nonoverlapping generations. Several papers (see [7, 8, 10, 12, 15], as well as Section 4 for an overview) have shown that the coalescent approximation applies more generally. Then the coalescent time turns into $N_e \sim N/c$ generations for some positive constant $c$ determined by the particular features of the population model under consideration. Following [15], we call $N_e$ the *coalescent effective size* (CES) of such a population model.

In the terminology of [17], the existence of a CES in the context of population genetic data is equivalent to a situation in which it is not possible to reject the basic WFM. This means that the coalescence pattern is indistinguishable from that of a WFM. The CES is a narrower version of the classical genetical concept of an *inbreeding effective size* (IES) designed to

Received July 2004; revised November 2004.

[1]Supported in part by The Bank of Sweden Tercentenary Foundation.

compare the rate of genetic drift in a given model with that of the WFM (see [5]). The existence of a CES implies that the IES exists as well and takes the same value. The reverse is not true: an IES may exist while a CES is absent (as in the case of convergence to a coalescent with multiple mergers [16]).

This paper looks for a general CES formula in the case of overlapping generations. The best known constant size genetic population model with overlapping generations is the Moran model, assuming that each unit of time one individual is killed and another produces an offspring, so that the population size $N$ remains constant over time. It is straightforward to verify that in this case the coalescent approximation holds with the coalescent time unit equal to $N^2/2$ units of the Moran model time. Since the generation time in the Moran model is $N$, this implies the existence of a CES with $N_e \sim N/2$. Here and elsewhere, $x_N \sim y_N$ means that $x_N/y_N \to 1$ as $N \to \infty$.

An age-structured version of the WFM, introduced in [4], discerns among $A$ age groups of constant sizes $N(a) = q(a)N, a = 1, 2, \ldots, A$. In contrast to the basic WFM, the full backward description of an age-structured ancestral process should include age-labelling of lineages. In terms of the probability $p(a)$ that a randomly chosen new-born individual has a parent of age $1 \leq a \leq A$, $p(1) + \cdots + p(A) = 1$, the lineage back of an individual in the age-structured WFM exhibits ages with the transition probabilities

$$(1) \quad \begin{bmatrix} p(1) & p(2) & p(3) & \cdots & p(A-1) & p(A) \\ 1 & 0 & 0 & \cdots & 0 & 0 \\ 0 & 1 & 0 & \cdots & 0 & 0 \\ \cdots & \cdots & \cdots & \cdots & \cdots & \cdots \\ 0 & 0 & 0 & \cdots & 1 & 0 \end{bmatrix},$$

the subdiagonal ones of course mirroring that an $(a+1)$-aged individual is viewed as stemming from an $a$-aged the preceding year, namely, herself. It is easy to verify that

$$\bar{\gamma} = (\gamma(1), \ldots, \gamma(A)),$$

with

$$(2) \quad \begin{aligned} \gamma(a) &= \frac{1}{\gamma}(p(a) + \cdots + p(A)), \\ \gamma &= \sum a p(a), \end{aligned}$$

is the corresponding stationary distribution of ages (in analogy with the stable age distribution of branching processes or deterministic age-structured population models, cf. [6]). Clearly, $\gamma(1) = 1/\gamma$.

According to [4], the IES of the age-structured WFM is $N_e \sim N/(c_{\text{age}} \gamma)$, where

$$(3) \qquad c_{\text{age}} = \sum_{a=1}^{A} \left( \frac{1}{q(a)} - \frac{1}{q(a-1)} \right) \gamma^2(a),$$

under the convention $\frac{1}{q(0)} = 0$. This IES formula takes into account that the generation time of the age-structured WFM is $\gamma$. In the particular case of constant fertility across ages, $p(a) \equiv q(a)$. Hence, with $\gamma(A+1) = 0$,

$$c_{\text{age}} = \sum_{a=1}^{A} \frac{\gamma^2(a) - \gamma^2(a+1)}{q(a)}$$

$$= \sum_{a=1}^{A} \frac{(q(a)/\gamma)(2\gamma(a) - (q(a)/\gamma))}{q(a)},$$

and the effective population size becomes $N_e \sim N/(2 - \gamma^{-1})$. This, in turn transforms to the formula for the Moran model as $\gamma \to \infty$, though the Moran model certainly has no fixed age distribution.

Just like the classical WFM, the offspring number of an $a$-aged individual in a large $(N \to \infty)$ age-structured WFM is asymptotically Poisson with the mean

$$(4) \qquad m(a) = p(a)q(1)/q(a).$$

In Section 2 we introduce an age-structured model allowing an arbitrary marginal reproduction law compatible with the constant population size assumption. The subject of our interest, the ancestral process of the age-structured population, is described at the end of Section 2. Our coalescent approximation result, Theorem 3.1, stated in Section 3, gives a CES formula for populations with exchangeable reproduction, which extends the CES formula $N_e \sim N/(c_{\text{age}} \gamma)$ for the age-structured WFM.

In Section 4 we interpret our CES formula in terms of earlier known formulae for geographically-structured WFMs with strong migration and for exchangeable populations with rapidly fluctuating sizes. Special attention is paid to the question whether $N_e$ is smaller than $N$. The final part of the paper is devoted to the proof of Theorem 3.1.

## 2. An age-structured population model.

Time is considered discrete with a unit to be called a year, for convenience. Let $a = 1, \ldots, A$ stand for the age of an individual, where $A < \infty$ is the maximal possible age. Each year the population has the same size $N$ and the age-composition also remains fixed,

$$\bar{N} = (N(1), \ldots, N(A)),$$

so that

$$N = N(1) + \cdots + N(A), \qquad N(1) \geq \cdots \geq N(A).$$

Individuals are assumed similar/exchangeable at least in the weak sense that all individuals of the same age have the same probabilities of surviving a year and have the same offspring number distribution. By the assumption of a fixed age structure, individuals can certainly not be independent of each other: if I survive, your chances diminish, and, similarly, if you have many kids a year, I tend to have few. But in some sense individuals should be interchangeable and we shall impose more of proper exchangeability, where needed. It is a good idea to visualize entities like the different individuals' lifespans or reproductions at various ages as exchangeable throughout, even though this is not always needed for the results.

The age structure determines the age distribution in the population through $q(a, N) := N(a)/N$, and even the life span distribution: the probability of surviving year $a$, given that you have survived the preceding year, is

$$\frac{N(a+1)}{N(a)} = \frac{q(a+1, N)}{q(a, N)}.$$

The survival function is the same for all individuals and determined by the products of yearly survival probabilities. If $L$ denotes individual life span, thus

$$\mathrm{P}(L \geq a) = \frac{q(a, N)}{q(1, N)}, \qquad a = 1, \ldots, A,$$

and

$$\mathrm{E}(L) = \frac{1}{q(1, N)}.$$

We assume that

(5)                    $$q(a, N) \to q(a) > 0, \qquad N \to \infty.$$

The vector of parameters $\bar{q} = (q(1), \ldots, q(A))$ then also describes the asymptotic life span distribution in large populations,

$$\mathrm{P}(L \geq a) = \frac{q(a)}{q(1)}, \qquad a = 1, \ldots, A,$$

and

$$\mathrm{E}(L) = \frac{1}{q(1)}.$$

Denote the one-year offspring numbers of the $a$-aged individuals by $\{\nu_l(a)\}_{l=1}^{N(a)}$ for $a = 1, \ldots, A$, assumed to be independent across age classes and exchangeable within them. Dropping the $\nu$-suffix for simplicity, we write $m(a, N) := \mathrm{E}(\nu(a))$, and require that

$$(6) \qquad \sum_{l=1}^{N(a)} \nu_l(a) = N(a)m(a, N), \qquad 1 \le a \le A,$$

so that

$$N(1) = \sum_{a=1}^{A} N(a)m(a, N).$$

Again, assume that there is convergence

$$(7) \qquad m(a, N) \to m(a), \qquad N \to \infty,$$

and, hence, that

$$q(1) = \sum_{a=1}^{A} q(a)m(a).$$

This means that the $p(a) := \frac{q(a)}{q(1)}m(a)$, $a = 1, 2, \ldots, A$, sum to one and give the asymptotic, so-called stable, distribution of age at childbearing in a critical population (cf. [6], Section 8.4), that is, one where the mean offspring number per individual equals one. Indeed, for fixed $N$, the expected yearly number of children of $a$-aged mothers is $N(a)m(a, N)$. Since the total number of children born in a year is $N(1)$, the distribution of age at childbearing is

$$\frac{N(a)m(a, N)}{N(1)} = \frac{q(a, N)}{q(1, N)}m(a, N)$$

and

$$\mathrm{E}(\nu(1) + \cdots + \nu(L)) = \mathrm{E}\left( \sum_{l=1}^{A} \sum_{a=1}^{l} \nu(a) \mathbb{1}_{\{L=l\}} \right)$$

$$= \sum_{a=1}^{A} \mathrm{E}(\nu(a) \mathbb{1}_{\{L \ge a\}})$$

$$= \sum_{a=1}^{A} m(a, N) \frac{q(a, N)}{q(1, N)} = 1,$$

since $m(a, N)$ is precisely the expected offspring number one year of a *surviving* $a$-aged individual. Of course, matters can not possibly stand otherwise when population size remains constant.

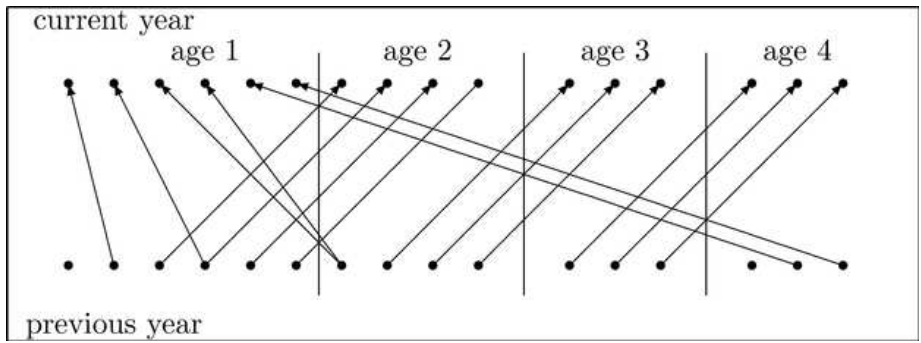

Fig. 1. *Forward picture.*

Having clarified the prospective view, we look backward, at the genealogy of $n$ individuals sampled out of the population the present year. Let $Z_0(a)$ denote the number of individuals of age $a$ among them, so that $Z_0 := Z_0(1) + \cdots + Z_0(A) = n$. The vector $\bar{Z}_0 = (Z_0(1), \ldots, Z_0(A))$ is the initial state of a Markov chain called the *ancestral process*. Its state at time $r$ is given by the numbers $\bar{Z}_r = (Z_r(1), \ldots, Z_r(A))$ of ancestors of the sampled individuals $r$ years ago, sorted by age. The total number of ancestors will be denoted by $Z_r := Z_r(1) + \cdots + Z_r(A)$. The Markov chain is time homogeneous and has a finite number of states. The class of states, where there is a single individual in one of the age classes, is absorbing.

Now, consider a situation where $k$ lines merge into one $a$-aged individual. If all the $k$ individuals are newborns, whose mother was of age $a$, we talk of a merger of type one, and denote it by $I(a, k)$. This situation can only occur if the mother is not in the sample. If she is, the $k$ merging lines will be those of $k - 1$ newborns and that from herself, one year later when she has attained age $a + 1$. Such mergers are said to be of type two and will be denoted $II(a, k)$.

EXAMPLE.    We illustrate the forward and backward views of this age-structured population model by two figures dealing with the case $N = 16$, $A = 4$, $N(1) = 6$, $N(2) = 4$, $N(3) = 3$, $N(4) = 3$. Figure 1 shows the development one year forward in time, the arrows indicating parent-offspring relationships and aging. The nonzero offspring sizes are $\nu_2(1) = \nu_4(1) = \nu_2(4) = \nu_3(4) = 1$; $\nu_1(2) = 2$.

Figure 2 presents the retrospective view of the model by tracing $n = 8$ ancestral lines one year back. We see one merger of type $I(2, 2)$ and one merger of type $II(1, 2)$. Each of the two mergers reduces the number of branches by one.

**3. Main result.** A first assumption on asymptotics, that age distributions $q(a, N)$ and expected offspring numbers $m(a, N)$ should converge as $N \to \infty$, has already been mentioned. In addition, we require that offspring variances stabilize,

$$(8) \qquad E(\nu^2(a)) \to m(a) + V(a), \qquad N \to \infty, \ 1 \le a \le A,$$

and that third moments do not grow too quickly,

$$(9) \qquad E(\nu^3(a)) = o(N), \qquad N \to \infty, \ 1 \le a \le A.$$

Here the limits $V(a) = \lim_{N \to \infty} E(\nu(a)(\nu(a) - 1))$ are never negative, since $\nu(a)(\nu(a) - 1) \ge 0$. [The case $V(a) = 0$ is not necessarily without interest, when there is an age-structure, since in this case individuals either may have given birth or not, $\nu(a) = 0$ or 1.] Whereas (8) serves to ensure that the time scale leading to the coalescent approximation is $T_N = N$, condition (9) is aimed at prohibiting multiple mergers of ancestral lines (cf. [16]).

THEOREM 3.1. *Assume* (5), (6), (7), (8), *and* (9)–(6), *interpreted to include the age-wise exchangebility and independence across ages mentioned before the equation itself.*

*Then the weak convergence to the Kingman coalescent*

$$(10) \qquad (Z_{\lfloor tN/\lambda \rfloor})_{t \ge 0} \to (R_t)_{t \ge 0}, \qquad N \to \infty$$

*holds with*

$$(11) \qquad \lambda = \frac{2}{\gamma} \sum_{a=1}^{A-1} \frac{p(a)\gamma(a+1)}{q(a)} + \frac{1}{\gamma^2 q^2(1)} \sum_{a=1}^{A} V(a) q(a)$$

*implying that a CES exists and satisfies*

$$(12) \qquad N_e \sim N/(\lambda\gamma), \ N \to \infty.$$

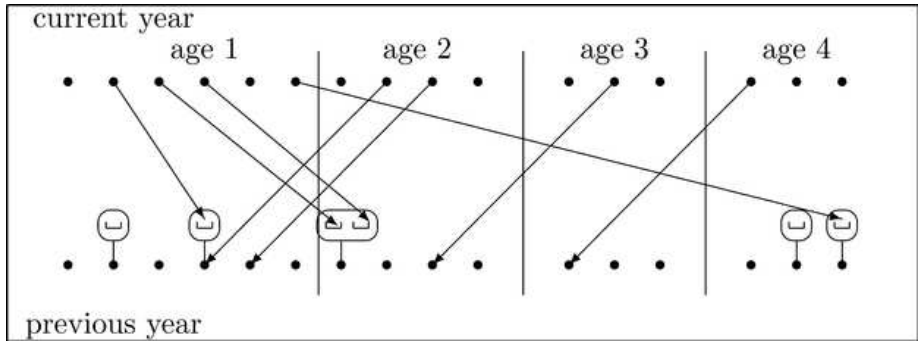

FIG. 2. *Backward picture.*

Our proof of (11) relies on asymptotics of joint factorial moments of off-spring numbers within an age class. It follows from Section 6 of [7] that the lower moments satisfy

$$(13) \qquad \mathrm{E}(\nu_1(a) \cdots \nu_j(a)) \to m^j(a)$$

and

$$(14) \qquad \mathrm{E}(\nu_1^{(2)}(a)\nu_2(a) \cdots \nu_j(a)) \to V(a)m^{j-1}(a),$$

while the higher factorial moments are bounded by

$$(15) \quad \mathrm{E}(\nu_1^{(x_1)}(a) \cdots \nu_j^{(x_j)}(a)) = o(N^{\delta-1}) \qquad \text{if } \delta := x_1 + \cdots + x_j - j \geq 2,$$

where $\nu^{(x)} := \nu(\nu-1) \cdots (\nu-x+1)$ denotes the descending factorial power. These relations amount to a sort of asymptotic independence of offspring numbers within age classes.

The result continues to hold if condition (6) is replaced by what is again a form of asymptotic independence, now between offspring numbers across age groups: for any $1 \leq m \leq A$, and any natural $j_1, \ldots, j_A$; $k_{11}, \ldots, k_{1j_1}; \ldots;$ $k_{A1}, \ldots, k_{Aj_A}$,

$$(16) \qquad \mathrm{E}\left(\prod_{a=1}^{A}\prod_{l=1}^{j_a} \nu_l^{k_{al}}(a)\right) \approx \prod_{a=1}^{A} \mathrm{E}\left(\prod_{l=1}^{j_a} \nu_l^{k_{al}}(a)\right), \qquad N \to \infty.$$

Here and elsewhere, $x_N \approx y_N$ means $x_N = y_N + o(1)$ and a product from one to zero equals one.

Though not necessary, the natural interpretation of the above setup is that, as $N \to \infty$, offspring numbers converge in distribution and in $L^p$ for any $p \geq 0$. The limiting random variables $\nu_l(a)$ will then satisfy (16) with equality. If they are bounded, it follows by the Weierstrass approximation theorem, that joint distribution functions factorize analogously. Hence, in the limit, reproduction of parents in different age classes is independent.

Such an extended version of Theorem 3.1 implies the existence of the CES of the age-structured WFM introduced in [4]. Indeed, the reproduction law of the age-structured WFM is given by the multinomial distribution

$$\mathrm{Mn}\left(N(1); \underbrace{\phi(1), \ldots, \phi(1)}_{N(1) \text{ times}}, \ldots, \underbrace{\phi(A), \ldots, \phi(A)}_{N(A) \text{ times}}\right),$$

where $\phi(a) = \frac{p(a)}{N(a)}$. Taking partial derivatives of the joint generation function

$$\mathrm{E}\left(\prod_{a=1}^{A}\prod_{l=1}^{N(a)} s_{al}^{\nu_l(a)}\right) = \left(\sum_{a=1}^{A}\sum_{l=1}^{N(a)} \phi(a)s_{al}\right)^{N(1)},$$

we obtain

$$\mathrm{E}\left(\prod_{a=1}^{A}\prod_{l=1}^{N(a)}\nu_l(a)^{(k_{al})}\right) \sim \prod_{a=1}^{A}\prod_{l=1}^{N(a)}(N(1)\phi(a))^{k_{al}} \sim \prod_{a=1}^{A}\prod_{l=1}^{N(a)}\mathrm{E}(\nu_l(a)^{(k_{al})}).$$

Thus, the offspring numbers are asymptotically independent both across and within different age classes, ensuring (16).

To see that $\lambda = c_{\mathrm{age}}$, as defined by (3), observe that exactly as in the argument following the definition,

$$
\begin{aligned}
c_{\mathrm{age}} &= \sum_{a=1}^{A}\frac{1}{q(a)}(\gamma^2(a) - \gamma^2(a+1)) \\
&= \frac{1}{\gamma}\sum_{a=1}^{A}\frac{p(a)}{q(a)}(\gamma(a) + \gamma(a+1)) \\
&= \frac{2}{\gamma}\sum_{a=1}^{A-1}\frac{p(a)\gamma(a+1)}{q(a)} + \frac{1}{\gamma^2}\sum_{a=1}^{A}\frac{p^2(a)}{q(a)}.
\end{aligned}
\tag{17}
$$

On the other hand, if the marginal distribution of the offspring number is asymptotically Poisson with mean (4), then $V(a) = (p(a)q(1)/q(a))^2$, turning the second term of (11) into that of the last expression above, so that $\lambda = c_{\mathrm{age}}$.

## 4. The coalescent rates $\lambda$, $c_{\mathrm{age}}$, $c_{\mathrm{geo}}$ and $c_{\mathrm{dem}}$.

The coalescent rate parameter $\lambda$ in (11) is a sum of two terms. From the derivation in Section 8, the first of these corresponds to a $\Pi(a, 2)$ merger, the second to one of type $I(a, 2)$. Notice that the second term disappears in the case $V(a) = 0$, with at most one offspring possible at the age $a$. In this section we interpret the two terms using the known CES formulae for the population models with 1. exchangeable reproduction, 2. strong migration and 3. fast size fluctuation.

The exchangeable haploid population model of [1, 2] is a flexible extension of the basic WFM, allowing an arbitrary marginal distribution of the offspring number $\nu$. According to [10], the CES of the exchangeable population satisfies $N_e \sim N/V_{\mathrm{hap}}$, where $V_{\mathrm{hap}} = \mathrm{E}(\nu(\nu-1))$ is the variance of the offspring number (the WFM corresponds to the symmetrical multinomial reproduction law with $V_{\mathrm{hap}} \approx 1$). For the haploid model, $V_{\mathrm{hap}}$ can be arbitrarily close to zero, so that no upper bound on $N_e$ is obtained. In [13] the last result was extended to diploid exchangeable models with random mating. If the haploid size of the diploid population is $N$, then it is shown that $N_e \sim 4N/V_{\mathrm{dip}}$, where again $V_{\mathrm{dip}} = \mathrm{E}(\nu(\nu-1))$ but now $\nu$ represents the number of diploid offspring to one couple. In this case $V_{\mathrm{dip}} \le 2$ and $N_e \le 2N$, with the upper bound reached when one couple produces exactly two children.

An important case of CES due to fast size fluctuations is considered in [7], where the demographic fluctuations backward in time occur according to a stationary Markov chain with the possible values $N(a) = q(a)N$, $a = 1, 2, \ldots, A$, the transition probabilities are $b_{ij}$, $i, j = 1, 2, \ldots, A$, and the stationary distribution is $(\gamma(1), \ldots, \gamma(A))$. For exchangeable reproduction, it is shown that $N_e \sim N/c_{\text{dem}}$ with

$$(18) \qquad c_{\text{dem}} = \sum_{i=1}^{A} \sum_{j=1}^{A} \gamma(i) b_{ij} V_{ij} q(j) q^{-2}(i),$$

where $V_{ij} = \mathrm{E}_{ij}(\nu(\nu - 1))$ measures variation in offspring numbers when the offspring generation size is $q(i)N$ and the parent generation size is $q(j)N$. This formula can be read as follows: two ancestral lines merge during a $q(i)N \to q(j)N$ backward size change at the rate equal to the rate $\gamma(i) b_{ij}$ of the size change times the conditional merger rate $V_{ij} q(j) q^{-2}(i)$. In the particular case of the WFM with fast fluctuations, the CES is approximated by the harmonic average of actual sizes $(\sum_{i=1}^{A} \frac{\gamma(i)}{q(i)N})^{-1}$, always smaller than the arithmetic average $N$ for nontrivial fluctuations.

Formula (18) yields the following interpretation of the second term in (11) corresponding to a $I(a, 2)$ merger. For two lines to merge at the age group $a$ as sister lines, they both should enter the age group $a$ immediately after visiting the age group 1 which happens at rate $(\gamma(1) b_{1a})^2 = (p(a)/\gamma)^2$. The corresponding conditional merger rate $V(a) q(a)/(m(a) q(a) N)^2$ is equal to that of a $m(a) q(a) N \to q(a) N$ backward size change. Multiplying these two rates and using (4) leads to the second term of (11).

To interpret the first term of (11) corresponding to a $II(a, 2)$ merger, we turn to the geographically structured WFM which is a key example in [15] illustrating the concept of CES. In this model a haploid population of constant size $N$ is split into $A$ subpopulations of constant sizes $N(a) = q(a)N$, so that $q(1) + \cdots + q(A) = 1$. Followed backward in time, an individual migrates from subpopulation $i$ to subpopulation $j$, with probability $b_{ij}$, and chooses its parent uniformly at random among $N(j)$ members of the parental subpopulation *independently* of other individuals. If migration is irreducible and aperiodic and $(\gamma(1), \ldots, \gamma(A))$ is the stationary distribution of the backward migration process, then, according to [15], the CES of the structured WFM satisfies $N_e \sim N/c_{\text{geo}}$, where

$$(19) \qquad c_{\text{geo}} = \sum_{a=1}^{A} \frac{1}{q(a)} \gamma^2(a).$$

As a formula for the inbreeding effective size (IES), (19) was discussed in [14], where it was pointed out that $c_{\text{geo}} \geq 1$, with the equality $c_{\text{geo}} = 1$ holding if only $\gamma(a) = q(a)$ for all $a$. The meaning of (19) is clear: for two lines to merge at the subpopulation $a$, they should be there at the same generation [rate $\gamma^2(a)$] and choose the same parent [rate $\frac{1}{Nq(a)}$].

The age-structured WFM is very similar to the geographically structured WFM with the transition probabilities $\|b_{ij}\|$ given by (1). However, (19) does not directly apply to the age-structured WFM, since individuals migrating backward in time from age group $a$ to age group $a-1$ sample their parents without replacement (thereby violating the assumption of independent choice of parents). Still, (19) helps in interpreting the first term of (11) corresponding to a $\Pi(a,2)$ merger. For two lines to merge at the age group $a$ as nonsister lines, they must enter the age group $a$ along different routes. One of the lines visits the age group $a$ immediately after visiting age group 1, which happens at rate $p(a)/\gamma$, while the other arrives through the age group $a+1$, which happens at rate $\gamma(a+1)$. [Notice that (19) also simplifies understanding of the last term in (17).]

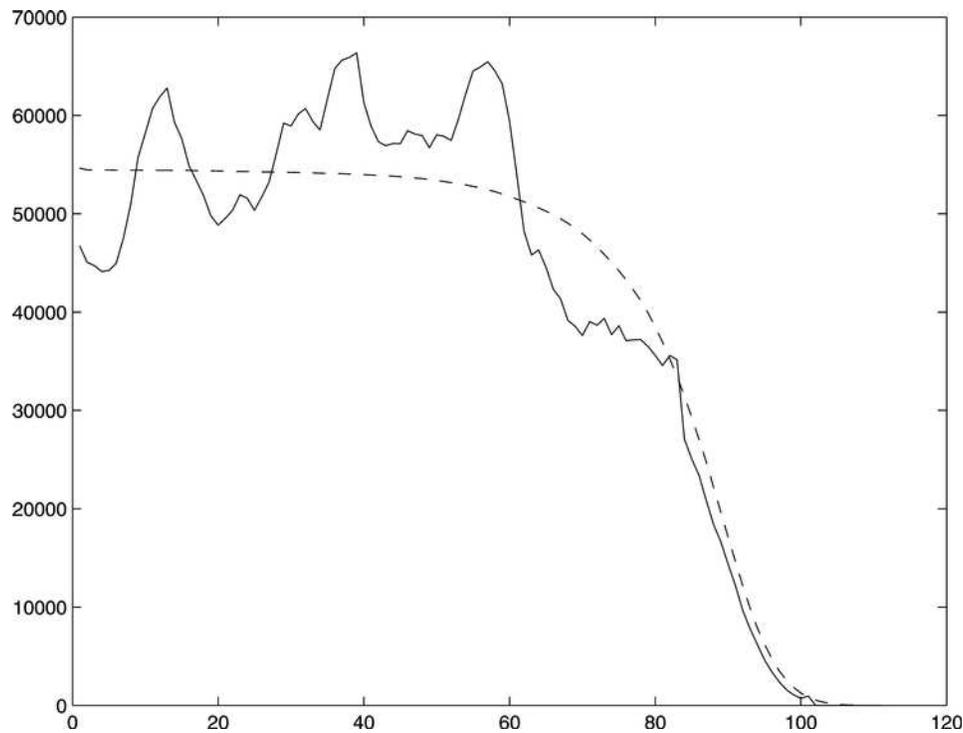

Fig. 3. *Observed (solid line) and stationary (dashed line) age group sizes of Swedish female population versus age.*

**5. Effective versus actual population size.** Recall that $q(1) \geq q(2) \geq \cdots \geq q(A)$ and $\gamma(1) = 1/\gamma$ to see that the first equality in (17) implies that

$$(20) \qquad \gamma c_{\text{age}} \geq \frac{\gamma}{q(1)} \sum_{a=1}^{A} (\gamma^2(a) - \gamma^2(a+1)) = \frac{1}{\gamma q(1)}.$$

Therefore, the CES of an age-structured WFM has the upper bound $N_e \leq \gamma q(1) N$, where the constant $\gamma q(1) = \gamma / \text{E}[L]$ is the ratio between the average age at child-bearing and the average life length. Similarly, since

$$(21) \qquad \frac{1}{\gamma} = \gamma \sum_{a=1}^{A} (\gamma^2(a) - \gamma^2(a+1)) = \sum_{a=1}^{A} p(a)(p(a) + 2\gamma(a+1)),$$

$$2 \sum_{a=1}^{A-1} p(a)\gamma(a+1) = \frac{1}{\gamma}\left(1 - \sum_{a=1}^{A} p^2(a)\right),$$

and we get a weaker upper bound

$$(22) \qquad N_e \leq \gamma q(1)\left(1 - \sum_{a=1}^{A} p^2(a)\right)^{-1} N$$

for the age-structured model with exchangeable reproduction.

These upper bounds could serve as fair estimates of CES in human and similar populations. For an illustration, we turn to Swedish official statistics for 2002 [18], yielding Figures 3 and 4, and $A = 111$, $\gamma = 30.6022$, $1/q(1) = 82.6094$.

The CES from (12) with $V(a) \equiv 0$ is $N_e = 0.3890N$. The age-structured WFM yields $N_e = 0.3677N$. (Cf. to Felsenstein's $N_e = 0.34N$ for U.S. population data, [4].) From Figure 4, $q(a) \approx q(1)$ for those $a$ where $\gamma(a)$ is not too small. This implies approximate equality in (20) and the above CES being close to the upper bounds (22), $N_e = 0.3919N$ and (20), $N_e = 0.3707N$.

**6. The transition probability.** In this section we derive an asymptotic formula for the one step transition probability

$$\Pi_{\bar{u}\bar{v}} := \text{P}(\bar{Z}_r = \bar{v} | \bar{Z}_{r-1} = \bar{u})$$

that a group of $u$ individuals with age distribution $\bar{u} = (u(1), \ldots, u(A))$ stems from a possibly smaller group of individuals from the previous year with age distribution $\bar{v} = (v(1), \ldots, v(A))$. We treat the parentage of newborns [$u(1)$ individuals of age 1] at time $r-1$ as balls to be allocated among $N$ boxes (potential mothers) at time $r$. Box $i$ in the age group $a$ contains a random number $\nu_i(a)$ of slots (see Figure 2), and each slot can accept one ball. The meaning of such an allocation is, of course, that one of the $u(1)$ individuals happens to be among the $\nu_i(a)$ children of the $i$th $a$-aged indivdual.

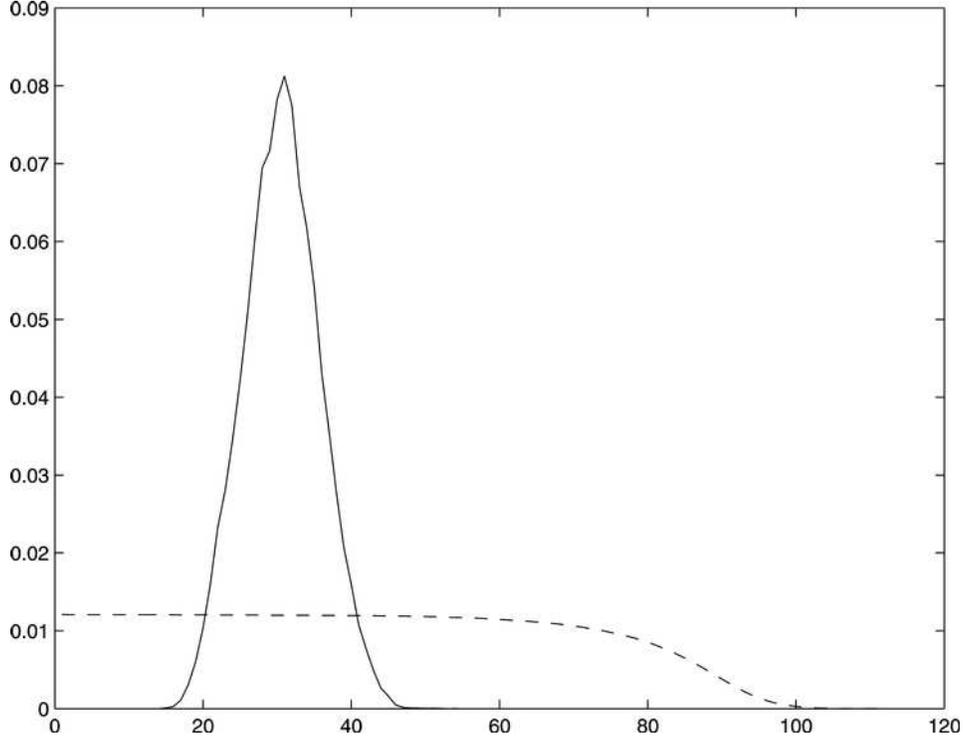

Fig. 4. *Solid line—$p(a)$, dashed line—$q(a)$, versus age $a$.*

Not newborn individuals, that is, of age $a + 1 \geq 2$ at time $r - 1$, just stem from themselves one year younger at time $r$, so $v(a) \geq u(a+1)$, writing $u(A+1) = 0$. Those remaining, $\alpha(a) = v(a) - u(a+1)$, must then have given birth to newborns. Combinatorically, we divide the $N(a)$ potential predecessors into $u(a+1)$ marked and $N(a) - u(a+1)$ unmarked boxes, in which balls can be placed to signify that the individual is among the predecessors. In Figure 2, for example, boxes are individuals in the previous year: counting from the left, we mark boxes 4 and 5 in the age group 1, box 3 in the age group 2 and, finally, box 1 in the age group 3.

Given the allocation result $\bar{v}$, we have $\alpha(a) = v(a) - u(a+1)$ unmarked boxes hosting at least one ball. We write

$$\bar{\alpha} = (\alpha(1), \ldots, \alpha(A)), \qquad u = \sum_{a=1}^{A} u(a), \qquad v = \sum_{a=1}^{A} v(a), \qquad \alpha = \sum_{a=1}^{A} \alpha(a),$$

$$\phi(\bar{\alpha}) = \binom{\alpha}{\alpha(1), \ldots, \alpha(A)} \prod_{a=1}^{A} p(a)^{\alpha(a)},$$

and notice that $\alpha = v - u + u(1)$.

The next step is to show that

$$\Pi_{\bar{u}\bar{v}} \sim (Nq(1))^{v-u} u(1)^{u-v} \phi(\bar{\alpha}) \left( \prod_{a=1}^{A} m(a)^{-\alpha(a)} \right)$$

$$(23)$$

$$\times \sum_{\bar{X}} \prod_{a=1}^{A} \mathrm{E}\left( \prod_{l=1}^{v(a)} \binom{\nu_l(a)}{x_l(a)} \right),$$

where $\bar{X} = (X(1), \ldots, X(A))$ and the vector $X(a) = (x_1(a), \ldots, x_{v(a)}(a))$ gives the numbers of balls in the $v(a)$ boxes. Numbers $x_l(a)$, with indices $1 \leq l \leq \alpha(a)$, correspond to unmarked boxes, while the indices $\alpha(a) + 1 \leq l \leq v(a)$ are meant for the marked boxes, so that summation in (23) is over all distinct arrays $\bar{X}$ satisfying

$$(24) \qquad \sum_{a=1}^{A} \sum_{l=1}^{v(a)} x_l(a) = u(1), \qquad x_l(a) \geq \mathbb{1}_{\{l \leq \alpha(a)\}}.$$

To illustrate the notation introduced in this section, we refer to Figure 2, which follows $u = 8$ ancestral lines one year back in the case $\bar{u} = (4, 2, 1, 1)$, $\bar{v} = (2, 2, 1, 1)$, $v = 6$. Here we have $\bar{\alpha} = (0, 1, 0, 1)$ and $X(1) = (1, 0)$, $X(2) = (2, 0)$, $X(3) = (0)$, $X(4) = (1)$.

PROOF OF (23). First we calculate the transition probability $\Pi_{\bar{u}\bar{v}}(\bar{I})$ when the set $\bar{I}$ of $v$ boxes is fixed. Here $\bar{I} = (I(1), \ldots, I(A))$ and the vector $I(a) = (i_1(a), \ldots, i_{v(a)}(a))$ lists the box positions taken from the set $\{1, \ldots, N(a)\}$. Again, positions $i_l(a)$, with indices $1 \leq l \leq \alpha(a)$, correspond to unmarked boxes, while the indices $\alpha(a) + 1 \leq l \leq v(a)$ are meant for the marked boxes. [In the case of Figure 2 we have $I(1) = (4, 5)$, $I(2) = (1, 3)$, $I(3) = (1)$, $I(4) = (3)$.] According to the allocation rules, we have

$$\Pi_{\bar{u}\bar{v}}(\bar{I}) = \frac{1}{\binom{N(1)}{u(1)} \prod_{a=1}^{A-1} \binom{N(a)}{u(a+1)}} \sum_{\bar{X}} \mathrm{E}\left( \prod_{a=1}^{A} \prod_{l=1}^{v(a)} \binom{\nu_{i_l(a)}(a)}{x_l(a)} \right)$$

$$\sim \frac{u(1)!}{N^{u(1)} q(1)^{u(1)} \prod_{a=1}^{A-1} \binom{N(a)}{u(a+1)}} \sum_{\bar{X}} \prod_{a=1}^{A} \mathrm{E}\left( \prod_{l=1}^{v(a)} \binom{\nu_l(a)}{x_l(a)} \right).$$

The expression becomes independent of $\bar{I}$ due to the exchangeability assumption. It remains to multiply the last expression with the number $\prod_{a=1}^{A} \binom{N(a)}{u(a+1)} \times \binom{N(a)-u(a+1)}{\alpha(a)}$ of ways to choose appropriate sets $\bar{I}$ of hosting boxes and do some simple algebra to obtain (23).

If $v \leq u - 2$, then (15) and (23) together imply that $\Pi_{\bar{u}\bar{v}} = o(N^{-1})$, meaning that multiple mergers of lineages are impossible in the limiting coalescent. □

**7. Transitions with $v \geq u - 1$.** We embark on a careful asymptotic analysis of the transition probability (23) with the most common transition type, when no ancestral lines merge: $v = u$. In this case $x_l(a) \equiv \mathbb{1}_{\{l \leq \alpha(a)\}}$, which in accordance with (13) entails $\Pi_{\bar{u}\bar{v}} \to A_{\bar{u}\bar{v}}$, where

$$(25) \qquad A_{\bar{u}\bar{v}} = \phi(\bar{\alpha}) \mathbb{1}_{\{\bar{u} \to \bar{v}\}} \mathbb{1}_{\{v = u\}}$$

and $\mathbb{1}_{\{\bar{u} \to \bar{v}\}} := \mathbb{1}_{\{v(1) \geq u(2), \ldots, v(A-1) \geq u(A)\}}$. The fact that the asymptotic transition probability takes the form of a multinomial distribution is easy to explain. As the set $\bar{u}$ of ancestral lines traced one year back does not change cardinality, the following happens. For $a \geq 2$, individuals of age $a$ turn to individuals of age $(a - 1)$, and lines from individuals of age 1 go to ages according to the multinomial Bernoulli scheme governed by the stationary distribution $\bar{p} = (p(1), \ldots, p(A))$ of individual's age at childbearing.

Let $\mathbf{A}_k = \|A_{\bar{u}\bar{v}}\|_{\bar{u}: u=k, \bar{v}: v=k}$ be the transition matrix at the level of $k$ ancestral lines. For $k = 1$, if the states are ordered as $(1, 0, \ldots, 0, 0), (0, 1, \ldots, 0, 0), \ldots,$ $(0, 0, \ldots, 0, 1)$, $\mathbf{A}_1$ is given by (1) with the stationary distribution (2). It is intuitively clear that, for an arbitrary $k$, the stationary distribution of the Markov chain with transition matrix $\mathbf{A}_k$ is given by the following multinomial distribution:

$$(26) \qquad \pi_k(\bar{u}) = \binom{k}{u(1), \ldots, u(A)} \prod_{a=1}^{A} \gamma(a)^{u(a)} \qquad \text{for } \bar{u} \text{ with } u = k.$$

Nevertheless, we give a formal proof of this fact using a computation technique which will be used later to produce less obvious results. Indeed, for any vector $\bar{v}$ with $v = k$, we have

$$\sum_{\bar{u}: u=k} \pi_k(\bar{u}) A_{\bar{u}\bar{v}} = \sum_{\bar{u}: u=k} \binom{k}{u(1), \ldots, u(A)} \prod_{a=2}^{A} \gamma(a)^{u(a)}$$
$$\times \binom{u(1)}{v(1) - u(2), \ldots, v(A)} \prod_{a=1}^{A} p(a)^{v(a) - u(a+1)} \mathbb{1}_{\{\bar{u} \to \bar{v}\}},$$

which equals

$$\binom{k}{v(1), \ldots, v(A)} \prod_{a=1}^{A} p(a)^{v(a)} \sum_{\bar{u}: u=k} \prod_{a=2}^{A} \binom{v(a-1)}{u(a)} \left(\frac{\gamma(a)}{p(a-1)}\right)^{u(a)} \mathbb{1}_{\{\bar{u} \to \bar{v}\}}.$$

Observe that the summation over $\bar{u}$ can be replaced by independent summations over the components $u(a)$, since

$$u(a) \leq v(a-1) \leq u(1) + u(a)$$
$$= u - u(2) - \cdots - u(a-1) - u(a+1) - \cdots - u(A)$$

and the summation index $u(a)$ is free to run between zero and

$$\min\{v(a-1), u-u(2)-\cdots-u(a-1)\} = v(a-1), \qquad a = 2,\ldots,A.$$

Therefore, the last sum converts to

$$\sum_{u(2)=0}^{v(1)}\cdots\sum_{u(A)=0}^{v(A-1)} \mathbb{1}_{\{u(1)=k-u(2)-\cdots-u(A)\}} \prod_{a=2}^{A} \binom{v(a-1)}{u(a)} \left(\frac{\gamma(a)}{p(a-1)}\right)^{u(a)}$$

and we conlude that

$$\sum_{\bar{u}:\,u=k} \pi_k(\bar{u}) A_{\bar{u}\bar{v}} = \binom{k}{v(1),\ldots,v(A)} \prod_{a=1}^{A} \gamma(a)^{v(a)} = \pi_k(\bar{v}).$$

Next we turn to the case $v = u - 1$ of exactly one pairwise merger and prove that $\Pi_{\bar{u}\bar{v}} \sim N^{-1} C_{\bar{u}\bar{v}}$, where

$$(27) \quad C_{\bar{u}\bar{v}} = \phi(\bar{\alpha}) \frac{u(1)}{q(1)} \sum_{a=1}^{A} \left[m(a)u(a+1) + \frac{\alpha(a)V(a)}{2m(a)}\right] \mathbb{1}_{\{\bar{u}\to\bar{v}\}} \mathbb{1}_{\{v=u-1\}}.$$

If $v = u - 1$, then $x_l(a) = \mathbb{1}_{\{l \le \alpha(b)\}} + \mathbb{1}_{\{a=b, l=j\}}$ for some $b \in \{1,\ldots,A\}$ and some $j \in \{1,\ldots,v(b)\}$, so that summation over $\bar{X}$ in (23) can be replaced by summation over $b$ and $j$. Note that indices $j \le \alpha(b)$ correspond to a single $I(b,2)$ case, while indices $j > \alpha(b)$ correspond to a single $II(b,2)$ merger. The corresponding components of (23) are computed with help of (13) and (14):

$$\lim_{N\to\infty} N\Pi_{\bar{u}\bar{v}}^{(1)} = \phi(\bar{\alpha}) \frac{u(1)}{q(1)} \prod_{a=1}^{A} m(a)^{-\alpha(a)}$$

$$\times \frac{1}{2} \sum_{b=1}^{A} \sum_{j=1}^{\alpha(b)} V(b) m(b)^{\alpha(b)-1} \prod_{a\neq b} m(a)^{\alpha(a)}$$

$$= \frac{1}{2} \phi(\bar{\alpha}) \frac{u(1)}{q(1)} \sum_{b=1}^{A} \frac{\alpha(b)}{m(b)} V(b)$$

and

$$\lim_{N\to\infty} N\Pi_{\bar{u}\bar{v}}^{(2)} = \phi(\bar{\alpha}) \frac{u(1)}{q(1)} \prod_{a=1}^{A} m(a)^{-\alpha(a)} \sum_{b=1}^{A-1} \sum_{j=1}^{u(b+1)} \prod_{a=1}^{A} m(a)^{\alpha(a)} m(b)$$

$$= \phi(\bar{\alpha}) \frac{u(1)}{q(1)} \sum_{b=1}^{A-1} m(b)u(b+1).$$

These two parts put together confirm (27).

**8. Proof of weak convergence toward the coalescent.** Theorem 3.1 can be established following the proof of weak convergence to the coalescent presented in Section 5 of [15]. The proof is based on Theorem 2.12 of [3] and the next lemma from [11].

LEMMA 8.1 (Möhle's lemma). *If* $\mathbf{A}$ *is a stochastic matrix such that* $\mathbf{P} = \lim_{k\to\infty} \mathbf{A}^k$ *exists, then*

$$\lim_{N\to\infty}\left(\mathbf{A} + \frac{1}{N}\mathbf{C} + o\left(\frac{1}{N}\right)\right)^{[Nt]} = \mathbf{P} - \mathbf{I} + e^{t\mathbf{G}},$$

*where* $\mathbf{I}$ *is the identity matrix and* $\mathbf{G} := \mathbf{PCP}$.

So far we have computed the transition matrix $\mathbf{\Pi} := \|\Pi_{\bar{u}\bar{v}}\|$ decomposition

$$\mathbf{\Pi} = \mathbf{A} + \frac{1}{N}\mathbf{C} + o\left(\frac{1}{N}\right),$$

with $\mathbf{A} := \|A_{\bar{u}\bar{v}}\|$, $\mathbf{C} := \|C_{\bar{u}\bar{v}}\|$. The only remaining calculation is to find the coalescence rate at level $k$ defined as

$$c_k := \sum_{\bar{u}:u=k} \pi_k(\bar{u})H(\bar{u}),$$

where $H(\bar{u})$ is the coalescence rate when the ancestor process is in configuration $\bar{u}$

$$H(\bar{u}) = \sum_{\bar{v}:v=u-1} C_{\bar{u}\bar{v}}.$$

According to (27) and (26),

$$c_k = \frac{1}{q(1)} \sum_{\bar{u}:u=k} \sum_{\bar{v}:v=k-1} \frac{u(1)!}{(v(1)-u(2))!\cdots v(A)!} \prod_{b=1}^{A} p(b)^{v(b)-u(b+1)}$$

$$\times \frac{k!}{u(1)!\cdots u(A)!} \prod_{b=1}^{A} \gamma(b)^{u(b)} \sum_{a=1}^{A}\left[m(a)u(a+1) + \frac{\alpha(a)V(a)}{2m(a)}\right]\mathbb{1}_{\{\bar{u}\to\bar{v}\}}.$$

After switching the order of summation over $\bar{v}$ and $\bar{u}$, and then regrouping the terms, we obtain

$$c_k = \frac{k}{q(1)} \sum_{\bar{v}:v=k-1} \binom{k-1}{v(1),\dots,v(A)} \prod_{b=1}^{A} p(b)^{v(b)}$$

$$\times \sum_{\bar{u}:u=k} \gamma^{-u(1)} \prod_{b=2}^{A} \binom{v(b-1)}{u(b)} \left(\frac{\gamma(b)}{p(b-1)}\right)^{u(b)}$$

$$\times \sum_{a=1}^{A}\left[m(a)u(a+1) + \frac{V(a)}{2m(a)}(v(a)-u(a+1))\right]\mathbb{1}_{\{\bar{u}\to\bar{v}\}}.$$

Opening the square brackets, we split the expression in two terms $c_k = c'_k + c''_k$. The first term equals (after putting the summation over $a$ in front of other sums and representing the summation over $\bar{u}$ as a multiple sum)

$$c'_k = \frac{k}{q(1)} \sum_{a=1}^{A-1} m(a) \sum_{\bar{v}:\, v=k-1} \binom{k-1}{v(1),\ldots,v(A)} \prod_{b=1}^{A} p(b)^{v(b)}$$

$$\times \sum_{u(2)=0}^{v(1)} \cdots \sum_{u(A)=0}^{v(A-1)} u(a+1) \prod_{b=2}^{A} \binom{v(b-1)}{u(b)} \left(\frac{\gamma(b)}{p(b-1)}\right)^{u(b)}$$

$$\times \gamma^{-u(1)} \mathbb{1}_{\{u(1)=k-u(2)-\cdots-u(A)\}}.$$

Since

$$(28) \qquad \sum_{i=0}^{n} \binom{n}{i} \left(\frac{\gamma(b)\gamma}{p(b-1)}\right)^i = \left(\frac{\gamma(b-1)\gamma}{p(b-1)}\right)^n,$$

we have

$$c'_k = \frac{k}{q(1)\gamma^k} \sum_{a=1}^{A-1} m(a) \sum_{\bar{v}:\, v=k-1} \binom{v}{v(1),\ldots,v(A)} \prod_{b\neq a} (\gamma(b)\gamma)^{v(b)} p(a)^{v(a)}$$

$$\times \sum_{u(a+1)=0}^{v(a)} u(a+1) \binom{v(a)}{u(a+1)} \left(\frac{\gamma(a+1)\gamma}{p(a)}\right)^{u(a+1)}.$$

Applying

$$\sum_{i=0}^{n} i \binom{n}{i} \left(\frac{\gamma(a+1)\gamma}{p(a)}\right)^i = n \left(\frac{\gamma(a+1)\gamma}{p(a)}\right) \left(\frac{\gamma(a)\gamma}{p(a)}\right)^{n-1}$$

to the last sum yields

$$c'_k = \frac{k}{q(1)\gamma} \sum_{a=1}^{A-1} m(a)\gamma(a+1)$$

$$\times \sum_{\bar{v}:\, v=k-1} \binom{k-1}{v(1),\ldots,v(A)} \prod_{b\neq a} \gamma(b)^{v(b)} v(a)\gamma(a)^{v(a)-1},$$

and the relation

$$(29) \qquad \sum_{i_1+\cdots+i_A=n} i_a \binom{n}{i_1,\ldots,i_A} \prod_{b\neq a} \gamma(b)^{i_b}\gamma(a)^{i_a-1}$$

$$= \frac{\partial}{\partial\gamma(a)} (\gamma(1)+\cdots+\gamma(A))^n = n$$

implies

$$c_k' = \binom{k}{2} \frac{2}{q(1)\gamma} \sum_{a=1}^{A-1} m(a)\gamma(a+1).$$

The second term is calculated similarly. From

$$c_k'' = \frac{k(k-1)}{2q(1)} \sum_{a=1}^{A} \frac{V(a)}{m(a)}$$

$$\times \sum_{\bar{v}\,:\,v=k-1} \binom{k-2}{v(1),\ldots,v(a-1),v(a)-1,v(a+1),\ldots,v(A)} \prod_{b=1}^{A} p(b)^{v(b)}$$

$$\times \sum_{u(2)=0}^{v(1)} \cdots \sum_{u(A)=0}^{v(A-1)} \prod_{b\neq a} \binom{v(b)}{u(b+1)} \binom{v(a)-1}{u(a+1)} \prod_{b=1}^{A-1} \left(\frac{\gamma(b+1)}{p(b)}\right)^{u(b+1)}$$

$$\times \gamma^{-u(1)} \mathbb{1}_{\{u(1)=k-u(2)-\cdots-u(A)\}}$$

and (28), we derive

$$c_k'' = \binom{k}{2} \frac{1}{q(1)\gamma^2} \sum_{a=1}^{A} \frac{V(a)}{m(a)}$$

$$\times \sum_{\bar{v}\,:\,v=k-1} \binom{k-2}{v(1),\ldots,v(a)-1,\ldots,v(A)} \prod_{b\neq a} \gamma(b)^{v(b)} \gamma(a)^{v(a)-1} p(a).$$

This together with (29) yields

$$c_k'' = \binom{k}{2} \frac{1}{q(1)\gamma^2} \sum_{a=1}^{A} \frac{V(a)}{m(a)} p(a).$$

We conclude that $c_k = \binom{k}{2}\lambda$, with

$$\lambda = \frac{2}{q(1)\gamma} \sum_{a=1}^{A-1} m(a)\gamma(a+1) + \frac{1}{q(1)\gamma^2} \sum_{a=1}^{A} \frac{V(a)}{m(a)} p(a).$$

In view of (4), this leads to (11).

**Acknowledgment.** We are grateful to the referee for close reading and useful suggestions.

## REFERENCES

[1] Cannings, C. (1974). The latent roots of certain Markov chains arising in genetics: A new approach, I. Haploid models. *Adv. in Appl. Probab.* **6** 260–290. MR343949

[2] CANNINGS, C. (1975). The latent roots of certain Markov chains arising in genetics:
    A new approach, II. Further haploid models. *Adv. in Appl. Probab.* **7** 264–282.
    MR371430
[3] ETHIER, S. N. and KURTZ, T. G. (1986). *Markov Processes: Characterisation and
    Convergence.* Wiley, New York. MR838085
[4] FELSENSTEIN, J. (1971). Inbreeding and variance effective numbers in populations
    with overlapping generations. *Genetics* **68** 581–597. MR403719
[5] HARTL, D. L. and CLARK, A. G. (1997). *Principles of Population Genetics*, 3rd ed.
    Sinauer Associates, Sunderland, MA.
[6] JAGERS, P. (1975). *Branching Processes with Biological Applications.* Wiley, Chich-
    ester. MR488341
[7] JAGERS, P. and SAGITOV, S. (2004). Convergence to the coalescent in populations
    of substantially varying size. *J. Appl. Probab.* **41** 368–378. MR2052578
[8] KAJ, I., KRONE, S. and LASCOUX, M. (2001). Coalescent theory for seed bank mod-
    els. *J. Appl. Probab.* **38** 285–301. MR1834743
[9] KINGMAN, J. F. C. (1982). On the genealogy of large populations. *J. Appl. Probab.*
    **19A** 27–43. MR633178
[10] KINGMAN, J. F. C. (1982). Exchangeability and the evolution of large populations.
    In *Exchangeability in Probability and Statistics* (G. Koch and F. Spizzichino,
    eds.) 97–112. North-Holland, Amsterdam. MR675968
[11] MÖHLE, M. (1998). A convergence theorem for Markov chains arising in popula-
    tion genetics and the coalescent with selfing. *Adv. in Appl. Probab.* **30** 493–512.
    MR1642850
[12] MÖHLE, M. (1998). Robustness results for the coalescent. *J. Appl. Probab.* **35** 438–
    447. MR1641829
[13] MÖHLE, M. and SAGITOV, S. (2003). Coalescent patterns in exchangeable diploid
    population models. *J. Math. Biol.* **47** 337–352. MR2024501
[14] NAGYLAKI, T. (1980). The strong-migration limit in geographically structured pop-
    ulations. *J. Math. Biol.* **9** 101–114. MR661421
[15] NORDBORG, M. and KRONE, S. (2002). Separation of time scales and convergence to
    the coalescent in structured populations. In *Modern Developments in Theoretical
    Population Genetics* (M. Slatkin and M. Veuille, eds.) 194–232. Oxford Univ.
    Press.
[16] SAGITOV, S. (2003). Convergence to the coalescent with simultaneous multiple merg-
    ers. *J. Appl. Probab.* **40** 839–854. MR2012671
[17] SJÖDIN, P., KAJ, I., KRONE, S., LASCOUX, M. and NORDBORG, M. (2004). On the
    meaning and existence of an effective population size. Preprint. Available at
    http://www. math.uu.se/~ikaj/papers.html.
[18] STATISTICS SWEDEN (Statistiska centralbyrån) (2003). Befolkningsstatistik
    2002, del 4. Födda och döda, civilståndsändringar m.m. Available at
    http://www.scb.se/templates/ publdb/publikation_2725.asp&plopnr=961.

DEPARTMENT OF MATHEMATICAL STATISTICS
CHALMERS UNIVERSITY OF TECHNOLOGY
AND GÖTEBORG UNIVERSITY
SE-412 96 GOTHENBURG
SWEDEN
E-MAIL: serik@math.chalmers.se
E-MAIL: jagers@math.chalmers.se